\documentstyle[11pt]{article}

\input amssym.def
\input amssym

\def \a{{\frak a}}
\def \Ad{{\rm Ad}}
\def \Aut{{\rm Aut}}
\def \b{{\frak b}}
\def \bs{{\backslash}}
\def \cd{{\rm cd}}
\def \C{{\Bbb C}}
\def \CE{{\cal E}}
\def \CF{{\cal F}}
\def \CH{{\cal H}}
\def \CO{{\cal O}}
\def \Eig{{\rm Eig}}
\def \Ell{{\rm Ell}}
\def \End{{\rm End}}

\def \g{{\frak g}}
\def \ga{\gamma}
\def \Ga{{\Gamma}}
\def \h{{\frak h}}
\def \Hom{{\rm Hom}}
\def \k{{\frak k}}
\def \la{{\lambda}}
\def \lap{\triangle}
\def \Lie{{\rm Lie}}
\def \m{{\frak m}}
\def \n{{\frak n}}
\def \N{{\Bbb N}}
\def \ord{{\rm ord}}
\def \p{{\frak p}}
\def \ph{\varphi}
\def \prf{{\bf Proof: }}
\def \P{{\Bbb P}}
\def \Per{{\rm Per}}
\def \qed{\hfill $\Box$

$ $

}
\def \Q{{\Bbb Q}}
\def \ra{\rightarrow}
\def \R{{\Bbb R}}
\def \supp{{\rm supp}}
\def \t{{\frak t}}
\def \tr{{\hspace{2pt}\rm tr\hspace{1pt}}}
\def \vol{{\hspace{1pt}\rm vol\hspace{1pt}}}
\def \Z{{\Bbb Z}}

\begin{document}

\title{Equivariant torsion of locally symmetric spaces} 
\author{Anton Deitmar\\ \small Math. Inst. d. Univ., INF 288, 69126 Heidelberg, Germany}
\date{}
\maketitle

\pagestyle{myheadings}
\markright{EQUIVARIANT TORSION...}

\tableofcontents

\newcommand{\rez}[1]{\frac{1}{#1}}
\newcommand{\der}[1]{\frac{\partial}{\partial #1}}
\newcommand{\binom}[2]{\left( \begin{array}{c}#1\\#2\end{array}\right)}

\newcounter{lemma}
\newcounter{corollary}
\newcounter{proposition}
\newcounter{theorem}

\newtheorem{conjecture}{\stepcounter{lemma} \stepcounter{corollary} 	
	\stepcounter{proposition}\stepcounter{theorem}Conjecture}[section]
\newtheorem{lemma}{\stepcounter{conjecture}\stepcounter{corollary}	
	\stepcounter{proposition}\stepcounter{theorem}Lemma}[section]
\newtheorem{corollary}{\stepcounter{conjecture}\stepcounter{lemma}
	\stepcounter{proposition}\stepcounter{theorem}Corollary}[section]
\newtheorem{proposition}{\stepcounter{conjecture}\stepcounter{lemma}
	\stepcounter{corollary}\stepcounter{theorem}Proposition}[section]
\newtheorem{theorem}{\stepcounter{conjecture} \stepcounter{lemma}
	\stepcounter{corollary}	\stepcounter{proposition}Theorem}[section]

$$ $$

\begin{center} {\bf Introduction} \end{center}

In search of an equivariant Riemann-Roch formula J.M Bismut \cite{Bism} and K. K\"ohler investigated in \cite{Koe} the equivariant torsion. The latter gave formulas expressing the equivariant torsion over $\P^n(\C)$ in terms of special values of zeta functions.
In this paper we will be concerned with equivariant torsion for Hermitian locally symmetric spaces.
It turns out that the equivariant torsion, or rather a quotient of the equivariant torsion over the equivariant $L^2$-torsion can be expressed as a special value of a zeta function which is defined in local geometric terms, i.e. by data related to closed geodesics and their Poincar\'e maps.
The $L^2$-contributions will vanish if the isometry has no fixed points.
Our formula may be considered as a Lefschetz formula in a wider sense.
More generally, consider an elliptic complex over a compact manifold.
The eigenvalue zero of the Laplacians gives the cohomology.
The nonzero eigenvalues give the torsion which may be viewed as the analogue of the Euler characteristic of the complex.
The analogue of the Atiyah $L^2$-index theorem \cite{At} would then be the equality of torsion and $L^2$-torsion. 
This however fails to hold. The quotient, which we express  by local data, measures the failure.

This paper deeply relates to \cite{Holtors} where we established such a formula in the non-equivariant setting, i.e. for the trivial isometry. 

To describe the results let $X$ be a Hermitian symmetric space of the noncompact type, $G$ the connected component of the isometry group of $X$ and $\Ga$ a cocompact neat discrete subgroup of the semisimple Lie group $G$.
Denote by $X_\Ga$ the compact quotient manifold $\Ga \bs X$.
Fix a holomorphic isometry $g_\Ga$ of $X_\Ga$ and a lift $g$ of $g_\Ga$ to $X$.
Suppose that $g$ lies in $G$.
The set of all lifts of $g_\Ga$ is given by the coset $\Ga g$ in $G$.
If $g_\Ga$ has fixed points the lift $g$ can be chosen of finite order.
We assume $g$ chosen in such a way.
The group $\Ga$ acts on the coset $\Ga g$ by conjugation.
Let $\Ell (\Ga g)$ denote the set of all elliptic $\Ga$-conjugacy classes in $\Ga g$, so a class $[\ga g]$ is in $\Ell (\ga g)$ if and only if $\ga g$ lies in a compact subgroup of $G$.
Note that $\Ell (\Ga g)$ is a finite set and is nonvoid if and only if $g_\Ga$ has fixed points.
Let $H$ be a Cartan subgroup of $G$ of splitrank $1$ and let $\CE_H(\Ga g)$ denote the set of nonelliptic $\Ga$-conjugacy classes in $\Ga g$, whose $G$-conjugacy classes meet $H$.
Let $\ph$ be a finite dimensional unitary representation of the group $\Ga$.
Since $\Ga$ is the fundamental group of $X_\Ga$ this defines a flat Hermitian vector bundle $E_\ph$ over $X_\Ga$.
For $\ga g\in\Ga g$ let $\Ga_{\ga g}$ denote the stabilizer in $\Ga$ under conjugation.
If $[\ga g]$ lies in $\CE_H(\Ga g)$ then the Euler characteristic of the group $\Ga_{\ga g}$ in the sense of \cite{Ser} will vanish.
There is however a notion of higher Euler characteristics the first of which is denoted $\chi_{_1}(\Ga_{\ga g})$ and this won't vanish in our case (see sec. \ref{sec_heat_trace}).
Let $A\subset H$ denote the (one dimensional) split torus and choose an ordering on the roots of $(G,A)$. 
The number of positive roots $c(H)$ will be $1$ or $2$.
Let $\n\subset \Lie G$ denote the sum of root spaces for the positive roots.
For $\Re (s)>>0$ consider the zeta function given in Weil form:
$$
Z_{H,1 ,\ph}^g(s) := \exp\left( -\sum_{[\ga g]\in \CE_H(\Ga g)}
	\frac{\chi_{_1}(\Ga_{\ga g}) \tr \ph(\ga)}
	     {\det (1-(\ga g)^{-1} | \n)}
	\frac{e^{-sl_{\ga g}}}
	     {\mu_{\ga g}}\right) ,
$$
where $\mu_{\ga g}$ is a certain multiplicity (see Theorem \ref{higher_heat_trace}).
Then we show that the function $Z_{H,1 ,\ph}^g$ extends to a logarithmic meromorphic function.
In the case that the order of the isometry $g_\Ga$ is $2$ we show that $Z_{H,1 ,\ph}^g$ extends to a meromorphic function and satisfies a Riemann hypothesis.
If furthermore $g_\Ga$ is fixed point free then $Z_{H,1 ,\ph}^g$ admits an Euler product expansion.

We will indicate in a special case how this zeta function can be interpreted as a geometric zeta function.
To this end assume that the isometry $g_\Ga$ has order 2 and no fixed points.
Then the group $\Ga'$ generated by $\Ga$ and $g$ is the fundamental group of the smooth quotient $X_{\Ga'}:=\Ga'\bs X$, of which $X_\Ga$ is a twofold cover.
Thus there is a bijection between the conjugacy classes in $\Ga'$ and the free homotopy classes of closed paths in $X_{\Ga'}$. Each class contains closed geodesics.
Lift those to $X$.
The space of geodesics on $X$ has a stratification by the rank.
The components of the lowest dimensional stratum correspond to the conjugacy classes of Cartan subgroups of splitrank one.
So the zeta function is defined by all closed geodesics on $X_{\Ga'}$ the lift of which lies in one component of the lowest dimensional stratum, and are such that they close in $X_{\Ga'}$ but not in $X_\Ga$.
The action on $\n$ can be interpreted as the expanding part of the Poincar\'e map and $\chi_{_1}(\Ga_{\ga g})$ is the orbifold Euler characteristic of the union of all closed geodesics in the class $[\ga g]$ modulo the action of the geodesic flow. 

Now we come to the equivariant torsion zeta function.
To keep the results neat we will assume that $G$ is simple, i.e. X is irreducible.
In the case that the number of roots $c(H)$ is $1$ let $Z_{H,\ph}^g :=Z_{H,1 ,\ph}^g(s+d(H)+b_0(H))$ for certain explicitly given constants $d(H), b_0(H)$.
There is a function $Z_{H,\ph}^g$ for the case $c(H)=2$ with similar properties.
Assume that $g_\Ga$ has order $2$ and define $Z_\ph^g$ by the finite product:
$$
Z_\ph^g(s) := \prod_{H/conjugation} Z_{H,\ph}^g(s),
$$
we obtain a meromorphic function on $\C$. We see that the order of $Z_\ph^g(s)$ at $s=0$ equals
$$
\sum_{q=0}^{\dim_\C X} q(-1)^q 
	\left( \tr (g_\Ga |\ker \lap_{0,q,\ph}) -
		 \dim \ph \sum_{[\ga g]\in \Ell (\Ga g)} 
			\tr_{\Ga_{\ga g}}(\ga g | \ker \tilde{\lap}_{0.q})
	\right),
$$
where $\lap_{p,q,\ph}$ is the Hodge Laplacian in $E_\ph$-valued $(p,q)$-forms and $\tilde{\lap}_{p,q}$ is the Hodge Laplacian on $X$.

We define $R_\ph^g(s) := Z_\ph^g(s)s^{-n_0}/c^g(X_\Ga ,\ph)$ for some explicitly given constant $c^g(X_\Ga ,\ph)$ and we obtain the special value
$$
R_\ph^g(0) = \frac{T_{g,hol}(X_\Ga ,\ph)}
		  {\prod_{[\ga g]\in \Ell(\Ga g)} T_{\ga g,hol}^{(2)}(X_\Ga)^{\dim \ph}},
$$
where $T_{g,hol}$ is the equivariant holomorphic torsion and $T_{\ga g,hol}^{(2)}$ the equivariant holomorphic $L^2$-torsion.

Concerning the restrictions we had to make we can say the following. 
If we drop the condition of $G$ being simple we will have to change the definition of the zeta functions drastically to obtain a similar theory.
Dropping the condition that $g_\Ga$ has order 2 results in $Z_\ph^g$ not being meromorphic anymore. 
Only the logarithmic derivative will be meromorphic then.
Nevertheless interpreting the special value as limit from above will yield the same results in that case, too.

\begin{center}{\bf Notation}\end{center}

We will write $\N ,\Z ,\Q ,\R ,\C$ for the sets of natural, integer, rational, real and complex numbers.

For any Hilbert space $\CH$ we write $B(\CH)$ for the algebra of all bounded linear operators on $\CH$.

For Lie groups $G,H,\dots$, we will write $\g_0,\h_0,\dots$ for the real Lie algebras and $\g ,\h ,\dots$ for their complexifications.

The convolution product on a Lie group $G$ will be denoted by $*$.
So if $f$ is a compactly supported smooth function on $G$ and $\ph$ is locally integrable then $\ph * f(x):= \int_G\ph(y)f(y^{-1}x)dy$, where $dy$ denotes a fixed Haar measure.

\section{Equivariant determinant and torsion}
\subsection{Classical}
Let $M$ denote a compact smooth Riemannian manifold and $E$ a smooth Hermitian vector bundle over $M$.
Let $D$ denote a {\bf generalized Laplacian} on $E$, i.e. $D$ is a second order differential operator with principal symbol $P_D(\xi)=|\xi |^2 Id$.
We further assume that $D$ is symmetric and semipositive.
(See \cite{BGV} for general information on such operators.)
Now let $g$ be an isometry of $M$ that lifts to a fibrewise linear isometry of $E$ denoted by the same letter.
Then $g$ acts on sections $s$ of $E$ by the pullback $g^*s(x) = g^{-1}(s(gx))$. 
Since $g$ is an isometry it acts unitarily  on the space of $L^2$-sections $L^2(E)$.
It leaves invariant the subspace of $C^\infty$-sections.
We assume that $g$ commutes with $D$.

We define the {\bf equivariant zeta function} of $D$ as
$$
\zeta_{g,D}(s) = \tr (g^*(D')^{-s})
$$
for $\Re (s) >>0$ and $D'=D|_{(\ker D)^\perp}$.
The existence of $\zeta_{g,D}$ is clear since some negative power of $D'$ is of trace class by Weyl's asymptotic law.

Note that for $g=1$ we have
$$
\zeta_{1,D}(s) =\zeta_D(s),
$$ 
the usual zeta function of the differential operator $D$.
Note further that, if $g$ acts on $M$ without fixed points and is of finite order, say $g^n = 1$ then 
$$
\rez{n} \sum_{k=1}^n \zeta_{g^k,D} = \zeta_{\tilde{D}},
$$ 
where $\tilde{D}$ is the pushdown of $D$ to the smooth quotient $<g>\bs M$.

We want to show that $\zeta_{g,D}$ extends to a meromorphic function on the entire plane.
To this end we consider the heat operator $e^{-tD}$.
It is known that $e^{-tD}$ is a smoothing operator and therefore has a smooth Schwartz kernel $<x|e^{-tD}|y>$, which is a section of the bundle $E \boxtimes E^*$, the exterior tensor product of $E$ with its dual over $M\times M$ with the property that $(E\boxtimes E^*)_{(x,y)} = \Hom_\C(E_x,E_y)$.
The operator $g^*e^{-tD}$ is also a smoothing operator with kernel $<x|g^*e^{-tD}|y> = g^{-1}<gx|e^{-tD}|y>$.

The fixed point set $M^g$ of $g$ is a submanifold of $M$. We have the {\bf small time asymptotics}:
$$
\tr g^*e^{-tD} \sim \sum_{k=0}^\infty c_k t^{k-\frac{\dim M^g}{2}},
$$
as $t\rightarrow 0$ for some coefficients $c_k$ which can be expressed as integral over $M^g$ of certain differential forms (see \cite{BGV}, p.193). 
Especially in the case $M^g=\emptyset$ it follows that $\tr\ g^*e^{-tD}$ is rapidly decreasing at $t=0$.
We are now able to prove

\begin{proposition}
The zeta function $\zeta_{g, D}$ extends to a meromorphic function on $\C$ with poles only at $s=(\dim\ M^g/2)-k$ for $k\geq 0$ an integer.
The residue of $\zeta_{g,D}$ at such a point is $c_k/\Ga (\frac{\dim\ M^g}{2}-k)$.
In particular, if $g$ is fixed point free then $\zeta_{g,D}$ is entire.
\end{proposition}
  
\prf
We have the identity 
$$
\zeta_{g,D}(s) = \rez{\Ga(s)}\int_0^\infty t^{s-1} \tr\ g^*e^{-tD'} dt
$$ 
for $\Re (s)>>0$ and $\tr\ g^*e^{-tD'} = \tr\ g^*e^{-tD} -\tr\ (g|\ker D)$.
Now split the integral as $\int_0^1 +\int_1^\infty$. The second summand converges for all $s$ and thus defines an entire function. In the first part substitute the asymptotic expansion to get the claim.
\qed

Note that by the proposition the zeta function $\zeta_{g,D}(s)$ is holomorphic at $s=0$.
We define the {\bf equivariant determinant} of $D$ as
$$
{\det}_g(D) := \exp(-\zeta_{g,D}'(0)).
$$
Note that by definition we have $\det_g(D) = {\prod_\lambda}^{reg} \la^{\tr (g|\Eig(\lambda ,D))}$, where the regularized product is taken over the nonzero eigenvalues of $D$.

For $\la > 0$ we now consider the operator $D+\la$. The above applies to this operator as well.

\begin{proposition}
The function $\la \mapsto \det_g(D+\la)$ extends to a holomorphic function on $\C -(-\infty ,0]$. We have $\det_g(D) = \lim_{\la \downarrow 0} \det_g(D+\la)\la^{-\tr (g|\ker D)}$.
\end{proposition}

\prf
The differential equation $\der{\la} \zeta_{g,D+\la}(s) = -s\zeta_{g,D+\la}(s+1)$ implies for $m\in \N$
$$
(\der{\la})^{m+1} \zeta_{g,D+\la}(s) = (-1)^{m+1} s(s+1) \dots (s+m) \zeta_{g,D+\la}(s+m+1),
$$
so that for $m$ large enough it follows $(\der{\la})^{m+1} \zeta_{g,D+\la}(0)=0$.

We have $\log \det_g(D+\la) = \lim_{s\ra 0}(\zeta_{g,D+\la}(s)-\zeta_{g,D+\la}(0))/s$ and the $s$-limit may be interchanged with the $\la$-derivation to give
$$
(\der{\la})^{m+1} \log {\det}_g(D+\la) = (-1)^{m+1} \Ga (m+1) \sum_{n=1}^\infty \frac{\tr (g|\Eig (\la_n,D))}{(\la +\la_n)^{m+1}}.
$$
From this the first claim follows. For the second replace $D$ by $D'$.
\qed

To define torsion suppose $E$ comes with a connection compatible with the Hermitian metric and let $\lap_{q,E}$ denote the Laplacian on $E$-valued $q$-forms.
The define the {\bf equivariant torsion} of $E$ as
$$
\tau_g (E) := \prod_{q=0}^{\dim M}{\det}_g(\lap_{q,E})^{q(-1)^q}.
$$

In this paper we will only be concerned with holomorphic equivariant torsion. 
So assume $M$ is K\"ahlerian, $E$ is a holomorphic bundle and $g$ is a holomorphic isometry.
Then $E$ has a unique connection compatible with the metric and holomorphic structure and we can define the Hodge-Laplace $\lap_{p,q,E}$ on $E$-valued $(p,q)$-forms.
The {\bf equivariant holomorphic torsion} is defined by
$$
T_{g,hol} (E) := \prod_{q=0}^{\dim_\C M} {\det}_g(\lap_{0,q,E})^{q(-1)^q}.
$$

\subsection{$L^2$-theory}
We will further be concerned with the $L^2$-versions of the above.
So consider now a compact smooth Riemannian manifold $X_\Ga$, its universal covering $X$ and its fundamental group $\Ga \subset \Aut (X)$.
Let $D_\Ga$ be a generalized Laplacian over a smooth Hermitian vector bundle $E_\Ga$ over $X_\Ga$, denote by $D$ and $E$ the corresponding lifts to $X$.
Now an isometry $g_\Ga$ of $X_\Ga$ will lift to an isometry of $X$ written $g$ such that $g\Ga g^{-1} =\Ga$.
Note that $g$ is only determined up to multiplication by elements of $\Ga$.

Choosing a fundamental domain $\CF$ of the $\Ga$-action on $X$ we get an isomorphism of $\Ga$-Hilbert modules
$$
L^2(E) \cong L^2(\Ga) \otimes L^2(E|_\CF) \cong L^2(\Ga)\otimes L^2(E_\Ga),
$$
where $\Ga$ acts on $\L^2(\Ga)$ by right translations and trivially on $L^2(E_\Ga)$. 
Let $VN(\Ga)\subset B(L^2(\Ga))$ denote the von Neumann algebra generated by the left translations $(L_\ga)_{\ga \in \Ga}$, then $VN(\Ga)$ coincides with the von Neumann algebra of all operators commuting with all right translations.
Therefore the algebra $B(L^2(E))^\Ga$ of all operators commuting with the $\Ga$-action becomes
$$
B(L^2(E))^\Ga \cong VN(\Ga) \otimes B(L^2(E_\Ga)).
$$
On the first tensor factor we have a canonical trace functional $\tau$ defined by $\tau(\sum_\ga c_\ga L_\ga) = c_e$ making $VN(\Ga)$ a type ${\rm II}_1$-von Neumann algebra if $\Ga$ is infinite.
Tensoring the trace $\tau$ with the usual trace on $B(L^2(E_\Ga))$ we get a trace $\tr_\Ga$ on $B(L^2(E))^\Ga$.
Let $T$ be an integral operator in $B(L^2(E))^\Ga$ with smooth kernel $<.|T|.>$ and a fundamental domain $\CF$ for $\Ga \bs X$, then a computation shows
$$
\tr_\Ga (T) = \int_\CF \tr <x|T|x> dx.
$$

The isometry $g$ acts on $L^2(E)$ but it is not $\Ga$-invariant.
Instead it is only invariant under the centralizer $\Ga_g$ of $g$ in $\Ga$.
But the $L^2$-theory works for $\Ga_g$ as well.
So as above we get a trace functional $\tr_{\Ga_g}$ on $B(L^2(E))^{\Ga_g}$.

\begin{lemma}
The operator $g^*e^{-tD}$ is of $\tr_{\Ga_g}$-trace class.
Its trace is
$$
\tr_{\Ga_g}(g^*e^{-tD}) = \int_{\CF_g}\tr g^{-1} <gx | e^{-tD} | x > dx,
$$
where $\CF_g$ is a fundamental domain for $\Ga_g \bs X$ and the integral converges absolutely.

This integral can also be written as the integral over the compact set $X_\Ga$ of the smooth function
$$
x \mapsto \sum_{\tau \in [g]_\Ga}\tr\tau^{-1} <\tau x|e^{-tD}|x>,
$$
where the sum runs over the $\Ga$-conjugacy class in $[g]_\Ga$ in $\Aut (X)$ of $g$.
\end{lemma}

\prf
The $\Ga$-invariance of $e^{-tD}$ amounts to $<\ga x | e^{-tD} |\ga y> = \ga <x | e^{-tD} |y >\ga^{-1}$, so it follows that the right hand side does not depend on the choice of $\CF_g$.
Now let $(u_i)_{i\in I}$ be a finite partition of unity on $X_{\Ga_g}=\Ga_g \bs X$, i.e. $u_i:X_{\Ga_g} \ra [0,1]$ is a smooth function for each i and $\sum_i u_i =1$.
The $u_i$ can be chosen such that for any $i,j\in I$ there is a fundamental domain $\CF_g$ such that no boundary point of $\CF_g$ maps to $\supp u_i \cup \supp u_j$. 
Let $\tilde{u}_i$ be the pullback of $u_i$ to $X$.
Since the $\tilde{u}_i$ give a partition of unity on $X$ it is sufficient to consider the operators $(g^*e^{-tD})_{i,j}$ with Schwartz kernel $<x|(g^*e^{-tD})_{i,j}|y>=\tilde{u}_i(x)\tilde{u}_j(y)<x|g^*e^{-tD}|y>$.
We will now make the isomorphism $L^2(E) \cong L^2(\Ga_g) \otimes L^2(E_{\Ga_g})$ more explicit.
So let $f\otimes \ph \in L^2(\Ga_g)\otimes L^2(E|_{\CF_g})$ then the corresponding element of $L^2(E)$ will be $\sum_{\ga \in \Ga_g}f(\ga)\ga^*\ph$.
The other way round we take $\psi \in L^2(E)$ and write it as $\psi = \sum_{\tau \in \Ga_g}\tau^*(\psi_\tau)$ with $\psi_\tau = {\bf 1}_{\CF_g}(\tau^{-1})^*\psi$.
With $T=(g^*e^{-tD})_{i,j}$ we write
\begin{eqnarray*}
T(f\otimes \ph) &=& \sum_{\ga \in \Ga_g}f(\ga)\ga^*T\ph\\
	&=& \sum_{\ga ,\tau}f(\ga)(\ga \tau)^*(T\ph)_\tau\\
	&=& \sum_{\tau\in \Ga_g}(R_{\tau^{-1}}f)\otimes(T\ph)_\tau.
\end{eqnarray*}

Therefore we conclude $\tr_{\Ga_g}(T)=\tr (\ph \mapsto (T\ph)_e)$, where $e$ is the neutral element of $\Ga_g$.
But the latter is the operator on $E|_{\CF_g}$ with kernel $<.|T|.>|_{\CF_g \times \CF_g}$.
By our assumption this gives a smooth kernel on $E_{\Ga_g}$ and the trace in question is just the trace of the integral operator thus defined.

We can build a fundamental domain $\CF_g$ out of $\CF$ by choosing a set of representatives $(\sigma)$ for $\Ga /\Ga_g$ and setting $\CF_g :=\bigcup_\sigma \sigma \CF$.
This induces an isomorphism $L^2(E_{\Ga_g}) \cong L^2(\Ga /\Ga_g) \otimes L^2(E_\Ga)$.
Using this isomorphism and the traces on the factors we see that everything boils down to showing that on $L^2(E|_\CF)$ the operator with kernel $\sum_{\sigma}<.|\sigma g e^{-tD}\sigma^{-1}|.>|_{\CF \times \CF}$ is of trace class and that its trace is the aboveclaimed.
This kernel can be written $\sum_{\tau \in [g]} \tau^{-1}<\tau x|e^{-tD}|y>$,
where the sum runs over the $\Ga$-conjugacy class of $g$.
Now growth estimates on the heat kernel as in \cite{CGT} show the absolute convergence of this sum locally uniformly in $x$ and $y$ and all its derivates.
Thus the kernel is smooth and since $X_\Ga$ is compact, the operator is of trace class and the trace is the integral over the diagonal, which gives the claim.
\qed

It is known that the small time asymptotics holds pointwise (compare \cite{BGV}) it follows that the trace $\tr_{\Ga_g}(g^*e^{-tD})$ also satisfies a similar asymptotics.
Unfortunately very little is known about large time asymptotics of $\tr_{\Ga_g}(g^*e^{-tD})$.
Let
$$
NS_g(D_\Ga) = \sup \{ \alpha \in \R | \tr_{\Ga_g} g^*e^{-tD'} = O(t^{-\alpha /2})\ {\rm as}\ t\ra \infty \}
$$
denote the {\bf equivariant Novikov-Shubin invariant} (compare \cite{GrSh}, \cite{LL}) of $D_\Ga$. 
J. Lott showed that the Novikov-Shubin invariants for $g=1$ of Laplacians are homotopy invariants.
J. Lott and W. L\"uck conjecture in \cite{LL} that the Novikov-Shubin invariants of Laplace operators are always positive rational or $\infty$.
In the situations we are going to consider in the sections to follow, where $X$ is a globally symmetric space of the noncompact type and $D$ is a motion invariant operator on $X$ Harish-Chandra's Fourier transform of orbital integrals shows that the equivariant Novikov-Shubin invariants are positive.

To proceed in the more general setting of the actual section we will have to {\bf assume} that $NS_g(D_\Ga)>0$.
We consider the integral
$$
\zeta_{g,D_\Ga}^1(s) := \rez{\Ga (s)} \int_0^1 t^{s-1} \tr_{\Ga_g} g^*e^{-tD'} dt,
$$
which converges for $\Re (s) >>0$ and extends to a meromorphic function on the entire plane which is regular at $s=0$ as follows from the small time asymptotics.
The integral
$$
\zeta_{g,D_\Ga}^2(s) := \rez{\Ga (s)} \int_1^\infty t^{s-1} \tr_{\Ga_g} g^*e^{-tD'} dt,
$$
converges for $\Re (s) < \rez{2} NS_g(D_\Ga)$. In this region we define the {\bf equivariant $L^2$-zeta function} of $D_\Ga$ as
$$
\zeta_{g,D_\Ga}^{(2)} (s) := \zeta_{g,D_\Ga}^1(s) +\zeta_{g,D_\Ga}^2(s).
$$

We define the {\bf equivariant $L^2$-determinant} of $D_\Ga$ as
$$
{\det}_g^{(2)}(D_\Ga) := \exp (-\frac{d}{ds}|_{s=0} \zeta_{g,D_\Ga}^{(2)}(s)).
$$

In analogy to the classical case we get

\begin{proposition}
The function $\la \mapsto \det_g^{(2)}(D_\Ga +\la)$, $\la >0$ extends to a holomorphic function on $\C -(-\infty ,0]$.
 
We have $\det_g^{(2)}(D_\Ga) = \lim_{\la \downarrow 0} \det_g^{(2)}(D_\Ga +\la)\la^{-\tr_{\Ga_g}(g|\ker D)}$.
\end{proposition}
\qed  

Again assume the bundle $E_\Ga$ comes with a compatible connection and let $\lap_{q,E_\Ga}$ be the Laplacian on $E_\Ga$-valued $q$-forms. The {\bf equivariant $L^2$-torsion} is by definition:
$$
\tau^{(2)}_g(E_\Ga) := \prod_{q=0}^{\dim X}{\det}_g^{(2)}(\lap_{q,E_\Ga})^{q(-1)^q}.
$$

Further, when $X$ again is K\"ahlerian, $E_\Ga$ holomorphic, $\lap_{p,q,E_\Ga}$ the Hodge-Laplacian on $E_\Ga$-valued $(p,q)$-forms then define
$$
T_{g,hol}^{(2)}(E_\Ga) := \prod_{q=0}^{\dim_\C X}{\det}_g^{(2)}(\lap_{0,q,E_\Ga})^{q(-1)^q}.
$$

\section{The equivariant trace of the heat kernel}
Now specialize to locally symmetric spaces.
Let $X$ be an Hermitian globally symmetric space without compact factors.
Write $G$ for the connected component of the group of isometries of $X$, then $G$ is a semisimple Lie group that acts transitively on $X$.
The stabilizer of any $x\in X$ is a maximal compact subgroup $K$ of $G$ so that by choosing a base point the space $X$ can be identified with $G/K$.
Let $\Ga \subset G$ be a discrete torsion free subgroup such that the quotient manifold $\Ga \bs G$ is compact.
Then $\Ga$ acts freely on $X$ and so $X_\Ga := \Ga \bs X$ is a compact Hermitian locally symmetric space.
Now let $g_\Ga$ be a nontrivial isometry of $X_\Ga$.
Then $g_\Ga$ lifts to an isometry $g$ of $X$ such that in the group of isometries of $X$ we have $g\Ga g^{-1} =\Ga$.
We will assume that $g$ belongs to the connected component $G$ of $Iso(X)$.
Note that only the coset $\Ga g$ in $G$ is determined by $g_\Ga$.

\begin{lemma} \label{finite_order}
There is a natural number $n$ such that $g^n\in \Ga$. In particular, the isometry $g_\Ga$ of $X_\Ga$ is of finite order.
\end{lemma}

\prf
Assume to the contrary that there is no $n\in \N$ with $g^n \in \Ga$. So the group $\Ga'$ generated by $\Ga$ and $g$ has infinite index over $\Ga$ hence $\Ga'$ can't be discrete.
Let $H$ be the closure of $\Ga'$ then $H$ is a Lie subgroup of $G$ of positive dimension and so the connected component $H^0$ is nontrivial.
Since $g$ normalizes $\Ga$ so does $\Ga'$ and hence $H$. 
Since $\Ga$ is discrete it follows that $H^0$ actually centralizes $\Ga$.
Now $G$ also has the structure of an algebraic group over $\R$ and in the present situation the group $\Ga$ is Zariski dense. Therefore $H^0$ also centralizes the entire group $G$ which is a contradiction since $G$ has trivial center.
\qed

Since $\Ga$ acts freely on the contractible space $X$ it follows that $X_\Ga$ is a $K(\Ga ,1)$-space, in particular, $\Ga$ is the fundamental group of $X_\Ga$.

The semisimple Lie group $G$ admits a compact Cartan subgroup $T\subset K$. 
We denote the real Lie algebras of $G,K,T$ by $\g_0,\k_0,\t_0$ and their complexifications by $\g ,\k ,\t$.
We will denote the Killing form of $\g$ by B.
As well, we will write $B$ for the diagonal of the Killing form, so $B(X) = B(X,X)$. Denote by $\p_0$ the orthocomplement of $\k_0$ in $\g_0$ with respect to $B$ then via the differential of \rm{exp} the space $\p_0$ is isomorphic to the real tangent space of $X=G/K$ at the point $eK$.
Let $\Phi (\t ,\g)$ denote the system of roots of $(\t ,\g)$, let $\Phi_c (\t ,\g) = \Phi (\t ,\k)$ denote the subset of compact roots and $\Phi_{nc} = \Phi -\Phi_c$ the set of noncompact roots.
To any root $\alpha$ let $\g_\alpha$ denote the corresponding root space.
Fix an ordering $\Phi^+$ on $ \Phi = \Phi (\t ,\g)$ and let $\p_\pm = \bigoplus_{\alpha \in \Phi_{nc}^+}\g_{\pm \alpha}$.
Then the complexification $\p$ of $\p_0$ splits as $\p = \p_+ \oplus \p_-$ and the ordering can be chosen such that this decomposition corresponds via the exponential to the decomposition of the complexified tangent space of $X$ into holomorphic and antiholomorphic part.

We will consider differential forms with coefficients in a flat Hermitian vector bundle. Such a bundle $E_\ph$ is given by a finite dimensional unitary representation $(\ph ,V_\ph)$ of the fundamental group $\Ga$.
More precisely $E_\ph = \Ga \bs X \times V_{\breve{\ph}}$ where $\Ga $ acts diagonally and $\breve{\ph}$ is the representation {\bf dual} to $\ph$.
We want the action of $g_\Ga$ on $X_\Ga$ to lift to $E_\ph$.
For this we {\bf assume} that $\ph^g =\ph$, where $\ph^g(\ga) = \ph(g\ga g^{-1})$.
The space of smooth $(p,q)$-forms with values in $E_\ph$ may be written as $\Omega^{p,q}(E_\ph) = (C^\infty (\Ga \bs G,\breve{\ph})\otimes \wedge^p\p_+ \otimes \wedge^q\p_-)^K$, where $C^\infty (\Ga \bs G,\breve{\ph}):= \{ s\in C^\infty(G,V_{\breve{\ph}})| s(\ga x) = \breve{\ph}(\ga)s(x) \}$.
In this setting the action of $g$ on $\Omega^{p,q}(E_\ph)$ is $g^*(s\otimes p_+ \otimes p_-) = g^*s \otimes p_+ \otimes p_-$ where $g^*s(x) = s(g x)$.

By \cite{BM} the heat operator $e^{-t\lap_{p,q,E_\ph}}$ has a smooth kernel $h_t^{p,q}$ of rapid decay in $(C^\infty(G)\otimes \End(\wedge^p\p_+ \otimes \wedge^q \p_-))^{K\times K}$.
For $t>0$ let
$$
f_t := \sum_{q=0}^{\dim_\C X}q(-1)^q \tr h_t^{0,q},
$$
where $\tr$ means the trace in $\End(\wedge^q\p_-)$.
Then $f_t$ acts by convolution from the left on 
$$
L^2(\Ga \bs G,\breve{\ph}):= \{s : G \ra V_{\breve{\ph}} | s(\ga x) = \breve{\ph}(\ga) s(x), \int_{\Ga \bs G} \parallel s(x) \parallel^2 dx <\infty\}
$$ 
and it follows 
$\tr(f_t|L^2(\Ga \bs G,\breve{\ph}))=\sum_{q=0}^{\dim_\C X}q(-1)^q \tr (e^{-t\lap_{0,q,E_\ph}})$.

\subsection{The equivariant trace formula}
Any smooth function $f$ of rapid decay on $G$ acts by convolution $s\mapsto s*f$ as a trace class operator on $L^2(\Ga \bs G,\breve{\ph})$.
The space $L^2(\Ga \bs G,\breve{\ph})$ is a unitary $G$-module and as such it decomposes into a discrete Hilbert sum
$$
L^2(\Ga \bs G,\breve{\ph}) = \bigoplus_{\pi \in \hat{G}}N_{\Ga ,\ph}(\pi)\pi,
$$
with finite multiplicities $N_{\Ga ,\ph}(\pi)$.

Therefore the trace of $f$ on the space $L^2(\Ga \bs G,\breve{\ph})$ is $\sum_{\pi \in \hat{G}} N_{\Ga ,\ph}(\pi) \tr \pi(f)$.
The Selberg trace formula asserts that this trace on the other hand equals
$$
\sum_{[\ga]} \vol(\Ga_\ga \bs G_\ga) \CO_\ga(f),
$$
 where the sum runs over all conjugacy classes $[\ga]$ in the group $\Ga$, the groups $\Ga_\ga$ and $G_\ga$ are the centralizers of $\ga \in \Ga$ in $\Ga$ and $G$ resp. and $\CO_\ga(f)$ is the {\bf orbital integral} $\int_{G_\ga \bs G} f(x^{-1}\ga x) dx$.

On the other hand, $g$ acts on $L^2(\Ga \bs G,\breve{\ph})$ by $g^*s(x):=s(gx)$ and this action commutes with $*f$ and the $G$-action from the right.
Therefore, for $\pi\in \hat{G}$ the isometry $g$ acts on $\Hom_G(\pi ,L^2(\Ga \bs G,\breve{\ph}))$ by a finite dimensional matrix $g_\pi^\Ga$.

\begin{proposition}
(Equivariant trace formula) Assume $\ph (g\ga g^{-1})=\ph (\ga)$ for $\ga \in \Ga$ then the trace of the operator $g^*(.*f)$ equals
$$
\sum_{\pi \in \hat{G}} \tr g_\pi^\Ga \tr \pi(f) = \sum_{[\ga g]_\Ga}
\vol(\Ga_{\ga g}\bs G_{\ga g}) \tr \ph(\ga)\CO_{\ga g}(f),
$$
where the sum runs over all $\Ga$-conjugacy classes in the coset $\Ga g$ and $\Ga_{\ga g}$ is the stabilizer of $\ga g$ in $\Ga$.
\end{proposition}

\prf
The left hand side clearly gives the trace of $g^*(.*f)$. To understand this operator, assume $s\in L^2(\Ga \bs G,\breve{\ph})$, then
\begin{eqnarray*}
g^*(s*f)(x)	&=& s*f(gx)\\
		&=& \int_G f(y^{-1} gx) s(y) dy\\
		&=& \sum_{\ga \in \Ga} \int_\CF f(y^{-1} \ga gx)\breve{\ph}(\ga^{-1})s(y) dy,
\end{eqnarray*}
where $\CF$ denotes a fundamental domain for the $\Ga$-action by left translates on $G$.
This shows that the operator $g^*(.*f)$ acts on $L^2(\Ga \bs G,\breve{\ph})$ as integral operator with smooth kernel
$$
k(x,y) = \sum_{\ga \in \Ga}f(y^{-1} \ga gx)\breve{\ph}(\ga^{-1}).
$$
Since $\Ga \bs G$ is compact, this is a trace class operator and its trace is given as the integral over the diagonal, i.e.:
\begin{eqnarray*}
\tr g(.*f) &=& \int_{\Ga \bs G} \tr k(x,x) dx\\
	&=& \sum_{\ga \in \Ga} \int_\CF f(x^{-1} \ga gx) \tr \ph (\ga)dx\\
	&=& \sum_{[\ga g]_\Ga} \int_{\Ga_{\ga g}\bs G} f(x^{-1} \ga gx) dx\tr \ph (\ga)\\
	&=& \sum_{[\ga g]_\Ga} \vol(\Ga_{\ga g}\bs G_{\ga g}) \tr \ph (\ga)\CO_{\ga g}(f).
\end{eqnarray*}
\qed

\subsection{The heat trace} \label{sec_heat_trace}
By Lemma \ref{finite_order} it follows that the group $\Ga'$ generated by $\Ga$ and $g$ again is a discrete subgroup of $G$.
It now may happen that the group $\Ga'$ admits torsion elements.
An element of $G$ is called {\bf elliptic} if it is contained in a compact subgroup.
Let $n\in \N$ denote the order of the isometry $g_\Ga$ of $X_\Ga$.

\begin{lemma}
The group $\Ga'$ has torsion elements if and only if some power $g_\Ga^k$ of $g_\Ga$ with $0<k<n$ acts on $X_\Ga$ with fixed points.

There are torsion elements in the coset $\Ga g \subset \Ga'$ if and only if $g_\Ga$ itself acts with fixed points. 
In that case the lift $g$ can be chosen of finite order.
\end{lemma}

\prf
Assume $\Ga'$ has torsion elements then there is a natural number $k<n$ and a $\ga \in \Ga$ such that $\ga g^k$ is torsion, hence elliptic, hence acts with fixed points on $X$ and so $g_\Ga^k$ acts with fixed points on $X_\Ga$.

The other way round if $g_\Ga^k$ has fixed points on $X_\Ga$ then there is $x\in X$ and $\ga \in \Ga$ such that $g^kx=\ga x$ so that $\ga^{-1}g^k$ is elliptic and since $\Ga'$ is discrete, $\ga^{-1}g^k$ is a torsion element.

The second assertion follows by the same proof restricted to the case $k=1$.
The fact that $g$ can be chosen of finite order is clear.
\qed

Any $\theta$-stable Cartan subgroup $H$ of $G$ will split as $H=AB$ where $A$ is the connected split component and $B$ is compact.
The dimension of $A$ is called the {\bf splitrank} of $H$.
Let $\a_0$ denote the real Lie algebra of $A$ and let $\a_0^+$ be the positive Weyl chamber as well as $A^+:= \exp (\a_0^+)$.
We will only be interested in Cartan subgroups of splitrank one.
For any such let $\CE_H(\Ga g)$ denote the set of $\Ga$-conjugacy classes $[\ga g]_\Ga$ such that $\ga g$ is in $G$ conjugate to an element $h=ab$ of $H$ with nontrivial split part $a$.
Let further $\Ell (\Ga g)$ denote the set of elliptic $\Ga$-conjugacy classes in $\Ga g$.
Note that $\Ell (\Ga g)$ is a finite set and is nonvoid if and only if $g_\Ga$ has fixed points on $X_\Ga$.

For a splitrank one Cartan $H=AB$ choose an order on the set of roots of $(\g ,\a)$ and let $c(H)$ the number of positive roots in the root system of $(\g ,\a)$, where $\a=\Lie_\C (A)$.
Then it is known (see \cite{Holtors}) that $c(H)$ equals $1$ or $2$.

In \cite{Ser} we find the notion of Euler characteristic $\chi(\Ga)$ for a lattice $\Ga$ in a reductive group $G$.
It has the property that $\chi(\Ga') = [\Ga :\Ga'] \chi(\Ga)$ for a sublattice $\Ga' \subset \Ga$.
Further, If $L$ has a nontrivial central split component then $\chi(\Ga)=0$ as is seen by fibre bundle arguments.

If $\Ga$ is torsion free it is of finite cohomological dimension.
Letting $b_j(\Ga)$ denote the Betti numbers we then define
$$
\chi_{_1}(\Ga) := \sum_{j=0}^{\cd(\Ga)} j(-1)^j b_j(\Ga).
$$
the {\bf first higher Euler number} of $\Ga$.

Recall that a lattice $\Ga$ is called {\bf neat} if for $\ga \in \Ga$ the adjoint $\Ad (\ga)\in \End(\g)$ does not have a root of unity as eigenvalue.
Every arithmetic $\Ga$ has a finite index neat subgroup \cite{Bor}.
 
We are going to need some notation from \cite{Holtors}.
At first, to our Cartan subgroup $H=AB$ there is a parabolic subgroup $P$ of $G$ such that $A$ is a split component of $P$, i.e. $P$ has Langlands decomposition $P=MAN$. Let $\m ,\a ,\n$ denote the corresponding Lie algebras over $\C$.

For an element $x$ of $G$ which is conjugate to some $A_x m_x \in AM$ we let $l_x$ denote the {\bf length} of $x$, i.e. $l_x :=\sqrt{B(\log a_x)}$.
This number is called length because for $\ga \in \Ga$ we have that $l_\ga$ is the length of any geodesic lying in the free homotopy class defined by $\ga$.

Let $L^M(\ga ,\tau)$ have the same meaning as in \cite{Holtors}.

For $l\geq 0$ let 
$$
b_l(H) :=(\frac{c(H)}{2} +\dim (\n) -1-l)|\frac{\alpha_r}{c(H)}|,
$$ 
where $\alpha_r$ is the real positive root in $\Phi (\g ,\h)$ and the absolute value comes from the Killing form. 
For $c(H)=2$ let 
$$
d_l(H):=\sqrt{b_l(H)^2+\frac{2B(\rho_{M_1,n})}{\dim(\p_{M_1,-})}+B(\rho_{M,n})}
$$
 where $M_1 :=M\cap G_1$ and $G=G_1\times G_2$ and $G_1$ is the simple factor of $G$ containing $A$.
In the case $c(H)=1$ we finally set 
$$
d(H) = \sqrt{B(\rho)-B(\rho_{K\cap G_2(H)})-B(\rho_{M_1(H)})}.
$$

The twisted trace formula together with the computation of the orbital integrals in \cite{Holtors} give

\begin{theorem} \label{higher_heat_trace}
Let $X_\Ga$ be a compact Hermitian locally symmetric space with neat fundamental group $\Ga$ and universal covering $X$ without compact factors.
Let $\ph$ denote a finite dimensional unitary representation of $\Ga$ and $E_\ph$ the flat vector bundle defined by its dual.
Write $\lap_{p,q,\ph}$ for the Hodge-Laplace operator on $E_\ph$-valued $(p,q)$-forms and $\lap_{p,q}$ for the Hodge Laplacian on $(p,q)$-forms on $X$. Let $g_\Ga$ be an isometry of $X_\Ga$ that lifts to an isometry $g$ in the connected component of the isometry group of $X$, then the trace of the operator $\sum_{q=0}^{\dim_\C X} q(-1)^q \tr g_\Ga^*e^{-t\lap_{0,q,\ph}}$ equals
\begin{eqnarray*} & \displaystyle
 \sum_{\begin{array}{c}H/conj.\\ c(H)=2\end{array}} \sum_{[\ga g] \in \CE_H(\Ga g)}
        \frac{\chi_{_1}(\Ga_{\ga g}) \tr \ph(\ga)}
                {\mu_{\ga g} \det(1-(\ga g)^{-1}|\n)}
        \frac{e^{-l_{(\ga g)}^2/4t}}{\sqrt{4\pi t}}
\\ & \displaystyle
        \sum_{l=0}^{\dim(\n_-)} (-1)^l e^{-td_l(H)^2}
        L^M(\ga g ,\wedge^l\n_-)
\\ & \displaystyle
+ \sum_{\begin{array}{c}H/conj.\\ c(H)=1\end{array}} \sum_{[\ga g] \in \CE_H(\Ga g)}
        \frac{\chi_{_1}(\Ga_{\ga g})\tr(\ph(\ga))}
             {\mu_{\ga g} \det(1-(\ga g)^{-1}|\n)}
 \frac{e^{-l_{(\ga g)}^2/4t}}{\sqrt{4\pi t}}
\\ & \displaystyle
                e^{-b_0(H)l_\ga} e^{td(H)^2} L^{G_2}(\ga g,1)
\\ & \displaystyle
        + \sum_{[\ga g]\in \Ell(\Ga g)}\tr \ph(\ga) \sum_{q=0}^{\dim_\C X} q(-1)^q \tr_{\Ga_{\ga g}} e^{-t\lap_{0,q}}.
\end{eqnarray*}
Here $\mu_{\ga g}:= l_{\ga g}/l_{(\ga g)_0^n}$ where $(\ga g)_0^n$ is the primitive element of $\Ga_{\ga g}$ underlying the element $(\ga g)^n$.
For unexplained notation we refer to \cite{Holtors}.
\end{theorem}

\prf
In the light of \cite{Holtors} there is only one thing that requires explanation. 
That is the occurrence of the Euler characteristic $\chi_{_1}(\Ga_{\ga g})$ and the factor $l_{(\ga g)_0^n}$. 
But these drop out of the volume factor in the trace formula using the arguments of \cite{Hitors}.
\qed

\section{Geometric zeta functions}
Now fix a $\theta$-stable splitrank one Cartan subgroup $H=AB$.
Further fix a parabolic subgroup $P$ with Langlands decomposition $P=MAN$.
Let $K_M:=K\cap M$ then $K_M$ is a maximal compact subgroup of the semisimple group $M$. The group $B$ then is a compact Cartan subgroup of $M$.
Let $\p_M$ be the positive part of the polar decomposition of the Lie algebra $\m$ of $M$.
Choose an order on the root system $\Phi (\m ,\b)$ where $\b := \Lie_\C B$ and let $\p_M = \p_{M,+}\oplus \p_{M,-}$ be the corresponding decomposition.

The symmetric space $X_M := M/K_M$ injects into $X$ and it inherits the holomorphic structure if and only if $c(H)=2$ (Lemma 2.2 in \cite{Holtors}). 
In that case we assume the order chosen such that $\p_{M,+}$ is mapped to the holomorphic tangent space.
Assuming $c(H)=2$ there is to every finite dimensional representation $(\tau ,V_\tau)$ of $K_M$ a compactly supported smooth function $g_\tau$ on $M$ such that for any $\xi \in \hat{M}$ we have
$$
\tr\xi(g_\tau) = \sum_{p=0}^{\dim_\C X_M}(-1)^p \dim (V_\xi \otimes \wedge^p\p_{M,-}\otimes V_{\breve{\tau}})^{K_M},
$$
where $\breve{\tau}$ is the dual to $\tau$ (see sec 2.4 in \cite{Holtors}).

Independent of the value of $c(H)$ there is to every finite dimensional representation $(\sigma ,V_\sigma)$ of the group $M$ a compactly supported smooth function $f_\sigma$ on $M$ such that for any $\xi \in \hat{M}$ we have
$$
\tr\xi(f_\sigma) = \sum_{p=0}^{\dim X_M} (-1)^p \dim (V_\xi \otimes \wedge^p \p \otimes V_{\breve{\sigma}})^{K_M}.
$$

We now recall

\begin{proposition}
Let $y$ be a semisimple element of the group $G$. If $y$ is not elliptic, the orbital integrals $\CO_y(f_\sigma)$ and $\CO_y(g_\tau)$ vanish. If $y$ is elliptic we may assume $y\in T$, where $T$ is a Cartan in $K$ and then we have
$$
\CO_y(f_\sigma) = \frac{{\tr\ \sigma(y)}|W(\t ,\g_y)| \prod_{\alpha \in \Phi_y^+}(\rho_y ,\alpha)}{[G_y:G_y^0]c_y},
$$
for all elliptic $y$ and
$$
\CO_y(g_\tau) = \frac{{\tr\ \tau(y)}}{\det(1-y^{-1} | \p_+)},
$$
if $y$ is regular elliptic. For general elliptic $y$ we have
$$
\CO_y(g_\tau) = \frac{\sum_{s\in W(T,K)} \det(s) \tilde{\omega}_y y^{s\la_{{\tau}}+\rho -\rho_K}}
        {[G_y :G_y^0]c_y y^\rho \prod_{\alpha \in \phi^+ - \phi_y^+}(1-y^{-\alpha})},
$$
where $c_y$ is Harish-Chandra's constant, it does only depend on the centralizer $G_y$ of $y$. Its value is given in \cite{Hitors}, further $\tilde{\omega}_y$ is the differential operator as in \cite{HC-DS} p.33.
\end{proposition}

\prf
\cite{Holtors}.
\qed

Let $G$ act on itself by conjugation, write $y.x = yxy^{-1}$, write $G.x$ for the orbit, so $G.x = \{ yxy^{-1} | y\in G \}$ as well as $G.S = \{ ysy^{-1} | s\in S , y\in G \}$ for any subset $S$ of $G$.
We are going to consider functions that are supported on the set $G.(MA)$. Let $\sigma$ be a finite dimensional representation of $M$.  
Now let
$l_a$ denote the length of $a$, so $l_a=|\log(a)|$.
Choose any $\eta : N \ra G$ which has compact support, is positive, invariant under $K\cap M$ and such that $\int_N \eta(n) dn =1.$ Let further $j$ be a natural number.

Given these data, let $\Phi =\Phi_{\eta, \sigma ,j,s} : G\ra \C$ be defined by
$$
\Phi(kn ma (kn)^{-1}) = \eta(n) f_\sigma(m) \frac{l_a^{j+1} e^{-sl_a}}{\det (1-(ma)^{-1}|\n)},
$$
for $k\in K$, $n\in N$, $m\in M$ and $a\in A$. Further $\phi(g)=0$ if $g$ is not in $G.(MA)$.

To see the welldefinedness of $\Phi$ recall first that by the decomposition $G=KP=KNMA$ every $g\in G.(MA)$ can be written as $kn ma (kn)^{-1}$.
Further suppose we have $kn ma (kn)^{-1}=k'n' m'a' (k'n')^{-1}$. Since $P$ is its own normalizer in $G$, it follows $k^{-1}k'\in K\cap P =K\cap M$. The invariance of $\eta$ and of $f_\sigma$ under $K\cap M$conjugation shows the welldefinedness.

We have the following generalization of Proposition 3.1 in \cite{Holtors}:

\begin{proposition}
The function $\Phi$ is $(j-\dim \n)$-times continuously differentiable.
For $j$ and $\Re (s)$ large enough it goes in the twisted trace formula and we have
$$
\sum_{\pi \in \hat{G}} N_{\Ga,\ph} (\pi) \tr g_\Ga^\pi \tr\pi(\Phi) = \sum_{[\ga g]_\Ga} \vol (\Ga_{\ga g} \bs G_{\ga g} ) \CO_{m_{\ga g}}^M(f_\sigma) \frac{l_{\ga g}^{j+1} e^{-sl_{\ga g}} \tr (\ph(\ga))}{\det (1-(\ga g)^{-1}|\n)},
$$
where the sum runs over all classes $[\ga g]_\Ga$ such that $\ga g$ is conjugate in $G$ to an element $m_{\ga g} a_{\ga g}$ of $MA^+$.
\end{proposition}
\qed

Besides the parabolic $P=MAN$ we also consider the opposite parabolic $\bar{P} =MA\bar{N}$. The Lie algebra of $\bar{N}$ is written $\bar{\n}$. Let $V$ denote a Harish-Chandra module of $G$ then we consider the Lie algebra homology $H_*(\bar{\n},V)$ and cohomology $H^*(\bar{\n},V)$. It is shown in \cite{HeSch} that these are Harish-Chandra modules of the group $MA$.

\begin{theorem} \label{L-funct1}
Let $X_\Ga$ be a compact Hermitian locally symmetric space with neat fundamental group $\Ga$ and universal covering $X$ without compact factors.
Let $\ph$ denote a finite dimensional unitary representation of $\Ga$.
Fix a splitrank $1$ Cartan $H=AB$ and a parabolic $P=MAN$.
Let $(\sigma ,V_\sigma)$ be a finite dimensional representation of $M$.
For $\Re (s)>>0$ define the function
$$
L_{H,\sigma ,\ph}^g(s) := \sum_{[\ga g] \in \CE_H(\Ga g)} \frac{\chi_{_1}(\Ga_{\ga g}) \tr \sigma (m_{\ga g}) \tr \ph(\ga) l_{(\ga g)_0^n}}{\det(1-(\ga g)^{-1}|\n)} e^{-sl_{\ga g}},
$$
where $n=\ord g_\Ga$ and $(\ga g)^n_0$ is the primitive element in $\Ga_{\ga g}$ underlying the element $(\ga g)^n \in \Ga_{\ga g}$.
Then $L_{H,\sigma ,\ph}^g$ extends to a meromorphic function with simple poles on the entire plane.
The residue at a point $s=\la (H_1)$, $\la \in \a^*$ is
$$
(-1)^{\dim \ \n} \sum_{\pi \in \hat{G}}N_{\Ga ,\ph}(\pi) \tr g_\Ga^\pi
\sum_{p,q}(-1)^{p+q} \dim (H^q(\bar{\n},\pi^0)\otimes \wedge^p\p_M \otimes V_{\breve{\tau}})^{K_M}_\la,
$$
where $(.)_\la$ denotes the generalized eigenspace and $H_1$ is the unique element of $\a_0^+$ with $B(H)=1$.
\end{theorem}

\prf
It follows from \cite{Hitors}
that for $\Ga$ neat
$$
\chi_{_1}(\Ga_{\ga g}) = 
\frac{\mid W(\g_{\ga g} ,\h)\mid \prod_{\alpha \in \Phi_{\ga g}^+}(\rho_{\ga g} ,\alpha )}
	{l_{(\ga g)_0^n} c_{\ga g} [G_{\ga g} :G_{\ga g}^0]} \vol (\Ga_{\ga g} \backslash G_{\ga g}),
$$
so that the geometric side of our trace formula will be
$$
\sum_{[\ga g]_\Ga\in \CE_H(\Ga g)} 
\frac{\chi_{_1}(\Ga_{\ga g}) \tr \sigma(m_{\ga g}) l_{(\ga g)_0^n} l_{\ga g}^{j+1} e^{-sl_{\ga g}}\tr(\ph(\ga))}
	{\det(1-(\ga g)^{-1} | \n)}
$$
which is just
$$
(-1)^{j+1} (\der{s})^{j+1} L_{H,\sigma,\ph}^g (s).
$$

The rest is analogous to the proof of Theorem 3.1 in \cite{Holtors}.
\qed

\begin{proposition}
The residue of $L_{H,\sigma ,\ph}^g(s)$ at $s=\la(H_1)$ is
$$
(-1)^{\dim N} \tr (g|\chi(\m \oplus \n ,K_M,C^\infty(\Ga \bs G ,\ph)\otimes \breve{\sigma})_{-\la}).
$$
\end{proposition}
\qed

For $\Re (s) >>0$ define
$$
Z_{H,\sigma ,\ph}^g(s) := \exp\left( -\sum_{[\ga g]\in \CE_H(\Ga g)}
	\frac{\chi_{_1}(\Ga_{\ga g}) \tr \sigma(m_{\ga g}) \tr \ph(\ga)}
	     {\det (1-(\ga g)^{-1} | \n)}
	\frac{e^{-sl_{\ga g}}}
	     {\mu_{\ga g}}\right) .
$$
Then the theorem asserts that the logarithmic derivative 
$$
L_{H,\sigma ,\ph}^g = (Z_{H,\sigma ,\ph}^g)'/Z_{H,\sigma ,\ph}^g
$$ 
extends to a meromorphic function on $\C$ with only simple poles.

An element $\ga g$ of $\Ga g$ is called {\bf $g$-primitive} if $\ga g = (\ga'g)^{kn+1}$ with $\ga'\in \Ga$ and $k\geq 0$ implies $k=0$.
If $g_\Ga$ has no fixed points then to every $\ga g$ there is a unique $g$-primitive $\Ga'g$ and $k\geq 0$ such that $\ga g = (\ga'g)^{kn+1}$.
Let $\CE_H^p(\Ga g)$ denote the subset of $\CE_H(\Ga g)$ formed by the $g$-primitive classes.

For any $y\in G$ let $\Per (\ga)$ denote the set of orders of nontrivial roots of unity occurring as eigenvalues of the adjoint $\Ad (y)\in \End(\g)$, so $m>1$ is in $\Per (y)$ if and only if a primitive $m$-th root of unity occurs as eigenvalue of $\Ad (y)$.
For any subset $I$ of $\Per (y)$ let $n_I$ denote the least common multiple of the elements of $I$.

\begin{theorem} \label{Z-funct1}
Let $Z_{H,\sigma ,\ph}^g$ be as above and now assume that $g_\Ga$ has order two.
Then $Z_{H,\sigma ,\ph}^g$ extends to a meromorphic function on $\C$.
The vanishing order of $Z_{H,\sigma ,\ph}^g(s)$ at $s=\la(H_1)$ is
$$
(-1)^{\dim \ \n} \sum_{\pi \in \hat{G}}N_{\Ga ,\ph}(\pi) \tr g_\Ga^\pi
\sum_{p,q}(-1)^{p+q} \dim (H^q(\bar{\n},\pi^0)\otimes \wedge^p\p_M \otimes V_{\breve{\tau}})^{K_M}_\la.
$$
If furthermore $g_\Ga$ acts on $X_\Ga$ without fixed points then $Z^g_{H,\sigma ,\ph}$ also admits an Euler product expansion:
$$
Z^g_{H,\sigma ,\ph}(s) = 
$$ $$
\prod_{[\ga g]\in \CE_H^p(\Ga g)} \prod_{N\geq 0}
			\prod_{I\subset \Per (\ga g)}
	\frac{\det (1-e^{-sl_{\ga g}} \ga \otimes \ga g | V_\ph \otimes V_N)^{a_I(\ga g)}}
	     {\det (1-e^{-2sl_{\ga g}} \ga^2 \otimes (\ga g)^2 | V_\ph \otimes V_N)^{a_I(\ga g)/2}},
$$
where $V_N := V_\sigma \otimes S^N(\n)$ and $\ga g$ acts on $V_N$ via $\sigma (m_{\ga g})\otimes \Ad^N(\ga g)$.

Further $a_I(\ga g) = \sum_{J\subset I}(-1)^{|I|+|J|}n_J \chi_{_1}(\Ga_{(\ga g)^{n_J}})/\mu_{(\ga g)^{n_J}}$.
\end{theorem}

In the case that besides the group $\Ga$ the group $\Ga'$ generated by $\Ga$ and $g$ also is neat we have $\Per (\ga g)=\emptyset$ and $a_\emptyset =\chi_{_1}(\Ga_{\ga g})/\mu_{\ga g}$.

\prf
If $g_\Ga$ has order two, so has $g_\Ga^\pi$ for any $\pi$ and thus the residues of $L_{H,\sigma,\ph}^g$ are integers.

The Euler product expansion is gotten as in \cite{Prodexp}.
\qed

Now assume $c(H)=2$ and fix a finite dimensional representation $(\tau ,V_\tau)$ of $K_M$.

\begin{theorem}
Let $\Ga$ be neat and $(\ph ,V_\ph)$ a finite dimensional unitary representation of $\Ga$.
Let $H$ be a $\theta$-stable Cartan subgroup of $G$ of splitrank $1$ with $c(H)=2$.
For $\Re (s) >>0$ define
$$
L_{H,\tau ,\ph}^{g,0}(s) := \sum_{[\ga g]\in \CE_H(\Ga g)} \frac{\chi_{_1}(\Ga_{\ga g})\tr \ph(\ga) L^M(\ga g ,\tau) l_{(\ga g)_0^n}}
	{\det(1-(\ga g)^{-1}|\n)} e^{-sl_{\ga g}}.
$$
Then $L_{H,\tau ,\ph}^{g,0}(s)$ extends to a meromorphic function with simple poles. The residue of $L_{H,\tau ,\ph}^{g,0}(s)$ at $s=\la(H_1)$ is
$$
-\sum_{\pi \in \hat{G}}N_{\Ga ,\ph}(\pi) \tr g_\Ga^\pi
\sum_{p,q}(-1)^{p+q} \dim (H^q(\bar{\n},\pi)\otimes \wedge^p\p_{M,-} \otimes V_{\breve{\tau}})_\la^{K_M},
$$
where $(.)_\la$ means the generalized $\la$-eigenspace.
\end{theorem}

\prf
The proof proceeds as the proof of theorem \ref{L-funct1} gut with $g_\tau$ taking the place of $f_\sigma$.
\qed

Again define
$$
Z_{H,\tau ,\ph}^{g,0}(s) := \exp \left( \sum_{[\ga g]\in \CE_H(\Ga g)} \frac{\chi_{_1}(\Ga_{\ga g})\tr \ph(\ga) L^M(\ga g ,\tau)}
	{\det(1-(\ga g)^{-1}|\n)} \frac{e^{-sl_{\ga g}}}{\mu_{\ga g}}\right) .
$$

As before we have $L_{H,\tau ,\ph}^{g,0} = (Z_{H,\tau ,\ph}^{g,0})'/Z_{H,\tau ,\ph}^{g,0}$.

\begin{theorem}
Assume that $g_\Ga$ has order two.
Then $Z_{H,\tau ,\ph}^{g,0}$ extends to a meromorphic function on $\C$. The vanishing order of $Z_{H,\tau ,\ph}^{g,0}$ at $s=\la (H_1)$ is
$$
-\sum_{\pi \in \hat{G}}N_{\Ga ,\ph}(\pi) \tr g_\Ga^\pi
\sum_{p,q}(-1)^{p+q} \dim (H^q(\bar{\n},\pi)\otimes \wedge^p\p_{M,-} \otimes V_{\breve{\tau}})_\la^{K_M}.
$$
\end{theorem}

\prf
As Theorem \ref{Z-funct1}.
\qed

Extend the definition of $Z_{H,\tau,\ph}^{g,0}(s)$ to arbitrary virtual representations of $K_M$ in the following way. 
Consider a finite dimensional virtual representation $\xi = \oplus_i a_i \tau_i$ with integers $a_i$ and $\tau_i\in \hat{K_M}$. 
Then let $Z_{H,\xi,\ph}^{g,0}(s) = \prod_i Z_{H,\tau_i,\ph}^{g,0}(s)^{a_i}.$

Now for $c(H)=2$ let 
$$
Z^g_{H,\ph}(s) = \prod_{l=0}^{\dim(\n_-)} Z^{g,0}_{H,\wedge^l\n_-,\ph}(s+d_l(H)+b_0(H))^{(-1)^l},
$$

In the case $c(H)=1$ let
$$
Z^g_{H,\ph}(s) = Z^g_{H,1,\ph}(s+d(H)+b_0(H)).
$$

\begin{proposition} \label{detformel}
Assume $G$ is simple, then for $\la >0$ we have the identity
$$
\prod_{q=0}^{\dim_\C X} \left( 
	\frac{\det_{g_\Ga} (\lap_{0,q,\ph}+\la)}
	     {\prod_{[\ga g]\in \Ell(\Ga g)}{\det}^{(2)}_{\ga g} (\lap_{0,q,\ph}+\la)}
\right)^{q(-1)^q}
            =
$$ $$
\prod_{\begin{array}{c}H/{\rm conj.}\\ \rm splitrank=1\\ c(H)=2\end{array}} \prod_{l=0}^{\dim \n_+}  \left( Z^g_{H,l,\ph} (b_0(H) + \sqrt{\la + d_l(H)^2}\right)^{(-1)^l}
$$ $$
\prod_{\begin{array}{c}H/{\rm conj.}\\ \rm splitrank=1\\ c(H)=1\end{array}}   Z^g_{H,1,\ph} (b_0(H) + \sqrt{\la + d(H)^2}).
$$
\end{proposition}

\prf
The equality is gotten by taking the Mellin transform of the expressions in Theorem \ref{higher_heat_trace}. \qed

Let $n_0$ be the order at $\la =0$ of the left hand side of the last proposition. Then the number $n_0$ equals
$$
 \sum_{q=0}^{\dim_\C X} q(-1)^q 
	\left( \tr (g_\Ga |\ker \lap_{0,q,\ph}) -
		 \dim \ph  \sum_{[\ga g]\in \Ell(\Ga g)}
			\tr_{\Ga_{\ga g}}(\ga g | \ker \tilde{\lap}_{0.q})
	\right).
$$
Now assume $g_\Ga$ has order $2$.
For $H$ a theta stable splitrank 1 Cartan with $c(H)=2$ and $l\geq 0$ let
$$
n^g_{H,l,\ph} = {\rm ord}_{s=b_0(H)+d_l(H)} Z^{g,0}_{H,\wedge^l\n_-,\ph}(s),
$$
so
$$
n^g_{H,l,\ph} = -\sum_{\pi \in \hat{G}}N_{\Ga ,\ph}(\pi) \tr g_\Ga^\pi
\sum_{p,q}(-1)^{p+q} \dim (H^q(\bar{\n},\pi)\otimes \wedge^p\p_{M,-} \otimes \wedge^l\n_+)_\la^{K_M}.
$$
Further, for $c(H)=1$ let $n_{H,\ph}^g$ be the order of $Z^g_{H,1,\ph}(s)$ at $s=b_0(H)+d(H)$, so
$$
n_{H,\ph}^g = -\sum_{\pi \in \hat{G}}N_{\Ga ,\ph}(\pi) \tr g_\Ga^\pi
\sum_{p,q}(-1)^{p+q} \dim (H^q(\bar{\n},\pi)\otimes \wedge^p\p_{M})_\la^{K_M}.
$$
 We then consider
$$
c^g(X_\Ga ,\ph)= (\prod_{H, c(H)=2}\prod_{l=0}^{\dim \n_+}(2d_l(H))^{n_{H,l,\ph}^g(-1)^l})
(\prod_{H, c(H)=1} (2d(H))^{n_{H,\ph}^g}).
$$

We assemble the results of this section to

\begin{theorem}
Assume that the isometry $g_\Ga$ has order two.
Then the zeta function $Z^g_{H,\ph}$ extends to a meromorphic function on the entire plane. Let
$$
Z^g_\ph(s) = \prod_{\begin{array}{c}H/{\rm conj.}\\ {\rm splitrank} = 1\end{array}} Z^g_{H,\ph}(s).
$$
Let $n_0$ be the order of $Z^g_\ph$ at zero then Proposition \ref{detformel} shows that the number $n_0$ equals
$$
\sum_{q=0}^{\dim_\C X} q(-1)^q 
	\left( \tr (g_\Ga |\ker \lap_{0,q,\ph}) -
		 \dim \ph  \sum_{[\ga g]\in \Ell(\Ga g)}
			\tr_{\Ga_{\ga g}}(\ga g | \ker \tilde{\lap}_{0.q})
	\right).
$$

Let $R^g_\ph(s) = Z^g_\ph(s)s^{-n_0}/c^g(X_\Ga ,\ph)$ then
$$
R^g_\ph(0) = \frac{T_{g,hol}(X_\Ga ,\ph)}{\prod_{[\ga g]\in\Ell(\Ga g)}T_{\ga g,hol}^{(2)}(X_\Ga)^{\dim\ph}}. 
$$ \qed
\end{theorem}

\tiny
Version: \today

\end{document}